\documentclass[reqno,12pt]{amsart}

\usepackage{amsmath}
\usepackage{amssymb}
\usepackage{latexsym}
\usepackage{color}

\allowdisplaybreaks[3]

\numberwithin{equation}{section}

\newtheorem{lemma}{Lemma}[section]

\newtheorem{prop}[lemma]{Proposition}
\newtheorem{thm}[lemma]{Theorem}
\newtheorem{cor}[lemma]{Corollary}

\newcommand{\Sptn}{ {\mathrm{Sp}}_{2n} }

\newcommand{\On}{{\mathrm O}_n}

\newcommand{\Gln}{\GL_n}
\newcommand{\Glnn}{\GL_n}
\newcommand{\Glk}{\GL_k }

\newcommand{\rM}{{\mathrm M}}
\newcommand{\C}{{\mathbb C}}
\newcommand{\R}{{\mathbb R}}
\newcommand{\Z}{{\mathbb Z}}
\newcommand{\Nn}{{U_{\mathrm{SO}_n}}}
\newcommand{\AS}{{A_{\mathrm{SO}_n}}}
\newcommand{\SO}{{\mathrm{SO}_n}}

\newcommand{\calb}{{\mathcal B}}
\newcommand{\calh}{{\mathcal H}}
\newcommand{\cale}{{\mathcal E}}

\newcommand{\calc}{{\mathcal C}}
\newcommand{\GL}{{ \mathrm{GL}} }
\newcommand{\calp}{{\mathcal{P}}}
\newcommand{\Mnk}{{\rM_{nk}}}
\newcommand{\pmnk}{{\calp(\Mnk)}}

\newcommand{\fraka}{{\mathfrak{A}}}

\newcommand{\fto}{{ \mathfrak{T}_{n,k,\ell} }}

\newcommand{\faol}{{ \fraka_{n,k,\ell} }}

\newcommand{\calhmkl}{{\calh(\rM_{n,k+\ell})}}

\newcommand{\LT}{{\mathrm{LM}}}
\newcommand{\calg}{{\mathcal{G}}}

\newcommand{\hs}{{\hspace{0.5in}}}
\newcommand{\nnn}{{\noindent}}

\newcommand{\Mnkl}{{\rM_{n,k+\ell}}}
\newcommand{\pmnkl}{{\calp(\Mnkl)}}

\newcommand{\ve}{{\varepsilon}}
\newcommand{\Mtnk}{{\rM_{2n,k}}}
\newcommand{\pmtnk}{{\calp(\Mtnk)}}
\newcommand{\asp}{{\mathcal{A}_{n,k,\ell}}}
\newcommand{\tG}{{\tilde{\Gamma}}}

\input xy
\xyoption{all}

\begin{document}

\title[Pieri Algebras]{Pieri algebras for the orthogonal and symplectic
groups}

\author {Sangjib Kim}

\address{
Department of Mathematics\\ The University of Arizona\\ 617 N. Santa
Rita Ave. P.O. Box 210089\\ Tucson, AZ 85721–0089 USA}
\email{sangjib@math.arizona.edu}

\author {Soo Teck Lee }

\thanks{{\tiny The second named author is  partially supported by NUS grant R-146-000-110-112.}}

\address {Department of Mathematics\\
National University of Singapore\\
2 Science Drive 2\\
Singapore 117543, Singapore.} \email {matleest@nus.edu.sg}

\begin{abstract}
We study the structure of a family of algebras which encodes an
iterated version of the  Pieri Rule for the complex orthogonal
group. In particular, we show that each of these algebras has a
standard monomial basis and has a flat deformation to a Hibi
algebra. There is also a parallel theory for the complex symplectic
group.
\end{abstract}


\maketitle

\section{Introduction}
  The standard monomial theory of Hodge, describing the homogeneous
coordinate ring of the Grassmannian and its natural extension to the case of the flag
manifold for $\Gln=\Gln(\C)$ (aka, the {\it flag algebra}),
is a landmark at the border between algebraic
geometry and representation theory. Hodge's results inspired
repeated attempts both to generalize them and to find an abstract
viewpoint from which they could be understood (\cite{Mu, La} and the
references therein). Important progress toward the second goal was
made by Gonciulea and Lakshmibai (\cite{GL}), who showed that the
$\Gln$ flag manifold could be described as a flat deformation of a
toric variety. Such a result has since been established for the flag
algebras of all reductive groups, and even for all multiplicity-free
actions (\cite{AB, Ca}).

\medskip
However, simply saying that a flag manifold is a flat deformation of
a toric variety does not capture all the structure inherent in
Hodge's original description, which also identified a finite family
of polynomial subrings of the $\Gln$ flag algebra, and showed that
the flag algebra was almost a direct sum of these subrings.

\medskip
An affine toric variety has a semigroup ring as its coordinate ring.
The semigroup in question has finite index in a {\it lattice cone} -
the intersection of ${\Z}^n$ with a rational polyhedral convex cone
in ${\R}^n$. In this context, it was  remarked in \cite{GL, Ki1, KM}
that the extra structure of Hodge's theory could be expressed by
describing the relevant cone (known as the Gelfand-Tsetlin cone, or equivalently
Gelfand-Tsetlin polytope) as a lattice cone attached to a Hibi ring:
the cone of non-negative integer-valued, order-preserving functions on a certain
partially ordered set (dubbed the {\it Gelfand-Tsetlin poset}
(\cite{Ho2, Ki1})).

\medskip
Hibi rings constitute an attractive class of algebras. Polynomial
rings are Hibi rings, and all Hibi rings share the concrete and
visible nature of polynomial rings. In particular, all Hibi rings
can be explicitly described in terms of generators and relations,
and they all have an explicit ``abstract standard monomial theory",
describing them as almost direct sums of certain polynomial subrings
(\cite{Ho2}). However, instead of a unique prototype in each
dimension, Hibi rings constitute a rich class with many examples in
each dimension.

\medskip It is therefore natural to ask, whether other rings besides
the $\Gln$ flag algebra occurring in representation theory can be
described as Hibi rings. A positive answer for the $\Sptn(\C)$ flag
algebra was already given in \cite{Ki1}. For the $\SO(\C)$ flag algebra,
a stable range case was studied in \cite{Ki2}.

\medskip
Standard monomial theory for $\Gln$ is closely linked with branching
- the restriction of representations of $\Gln$ to $\GL_{n-1}$, then
to $\GL_{n-2}$, and so forth.  A reciprocity law links branching
from $\GL_n$ to $\GL_{n-1}$ to decomposing a tensor product $\rho
\otimes S^m$, where $\rho$ indicates any irreducible representation
of $\GL_n$, and $S^m=S^m(\C^n)$ is a symmetric power of the standard
representation on ${\C}^n$. %
  The decomposition of such tensor products turns out to be
the same as the case of the Schubert calculus named after Mario Pieri,
so the result is known as the {\it Pieri rule}.

\medskip
Successive branchings correspond to decomposition of multiple tensor
products \begin{equation}\label{eq1}
\rho \otimes S^{m_1} \otimes S^{m_2} \otimes S^{m_3}
\otimes \cdots   \otimes S^{m_k}.
\end{equation}
Because of the connection with the Pieri rule, we call the algebra that
describes these tensor products   an {\it iterated Pieri algebra}, or simply
a {\it Pieri algebra}. The same techniques that show
that the flag algebra of $\GL_n$ is a flat deformation of a Hibi ring also establish that
a Pieri algebra is a flat deformation of a Hibi ring.

\medskip
In this paper, we consider the problem of describing certain
families of tensor products for the orthogonal and symplectic
groups, analogous to the products \eqref{eq1} for $\GL_n$. In the
manner of \cite{HL1, HTW3}, we construct an algebra that describes
these tensor products, which we again call a {\em Pieri algebra}.
Our main result is that each of these algebras has a standard
monomial basis and has a flat deformation to a  Hibi ring.

\medskip Here is a slightly more detailed overview of the paper.
The irreducible rational representations of $\On$ are labeled by
Young diagrams such that the sum of the lengths of the first two
columns is at most $n$ (\cite{Wy, GW, Ho1}). For such a
Young diagram $D$, let $\sigma^D_n$ be the corresponding irreducible
representation of $\On$. For positive integers $k$ and $\ell$ such
that $2(k+\ell)<n$ (we call this the stable range condition), we construct an algebra $\faol$ with the
following properties: $\faol$
 carries a multigrading, and  each of its homogeneous components can be
identified with the space of $\SO$ highest weight vectors of a
certain weight in a tensor product of the form
\begin{equation}\label{tp}
\sigma^D_n\otimes\sigma^{(p_1)}_n\otimes\sigma^{(p_2)}_n\otimes\cdots
\otimes\sigma^{(p_\ell)}_n \end{equation} where $D$ is a Young
diagram with at most $k$ rows and $p_1,...,p_\ell \geq 0$.   Thus,
the multiplicity of an irreducible representation of $\On$ in this
tensor product coincides with the dimension of the corresponding
homogeneous component. In this sense, the algebra encodes an
iterated version of the Pieri rule for $\On$, so  we   call the
algebra $\faol$ a  {\em Pieri algebra} for $\On$.

\medskip The structure of $\faol$ is closely related to a finite poset
$\tG=\tG(k,\ell)$ which arises from a consideration of the
multiplicities in tensor products of the form  \eqref{tp}. Let
$\Omega=\Omega(k,\ell)$ be the set of order preserving functions on
$\tG$ with nonnegative integral values. Then $\Omega$ can be
identified with a lattice cone. In particular, $\Omega$ is a
semigroup with respect to the usual addition of functions and it has
a distinguished finite set $\calg=\calg(k,\ell)$ of generators.
Using $\calg$, we identify a set $S=S(k,\ell)$ of generators for the
algebra $\faol$.
 The generating set
$\calg$ of $\Omega$ has a structure of a distributive lattice,  so
 it gives rise to a Hibi algebra (\cite{Hi, Ho2}).
This Hibi algebra is isomorphic to the semigroup algebra
$\C[\Omega]$ on $\Omega$. The partial ordering on $\calg$ also
induces a partial ordering on $S$.

\bigskip\noindent {\bf MAIN THEOREM.} {\em Let $n$, $k$ and
$\ell$ be positive integers such that $2(k+\ell)<n$.
\begin{enumerate}
\item[(a)] The algebra $\faol$ has a standard monomial theory for
$S$, that is, the  monomials in elements of the maximal chains in
$S$ form a basis for $\faol$.

\medskip
\item[(b)] There is a flat family of complex algebras with general fibre
$\faol$ and special fibre $\C[\Omega]$.
\end{enumerate}}

\bigskip\nnn We have also obtained the following multiplicity formula
from the structure of the algebras   $\faol$ :

\bigskip\nnn{\bf The Iterated Pieri Rule for $\On$.} {\em Let $k,\ell$ and $n$ be positive
integers such that $2(k+\ell)<n$, $D$ a Young diagram with at most
$k$ rows, and $P=(p_1,...,p_\ell)\in\Z^\ell_{\geq 0}$. For any Young
diagram $F$ which labels a representation of $\On$, the multiplicity
$m(F,D,P)$ of $\sigma^F_n$ in the tensor product
\[\sigma^D_n\otimes\sigma^{(p_1)}_n
\otimes\sigma^{(p_2)}_n\otimes\cdots \otimes\sigma^{(p_\ell)}_n\]
    is given by
\[m(F,D,P)= \sum_{E,A,B,C} K_{F/E,A} K_{D/E,B}\]
where  $K_{F/E,A}$ and  $K_{D/E,B}$ are skew Kostka numbers
(\cite{Sta}) and the sum is taken over all Young diagrams $E$ with
at most $k$ rows, $A=(a_i),B=(b_j)\in\Z^\ell_{\geq 0}$ and
$C=(c_{s,t})_{1\leq s<t\leq\ell}\in\Z^{\ell(\ell-1)/2}_{\geq 0}$
such that
\[p_i=a_i+b_i+c_{1i}+\cdots+c_{i-1,i}+c_{i,i+1}+\cdots+c_{i,\ell}\ \
 (1\leq i\leq \ell).\]}

\medskip Alternatively, the multiplicity $m(F,D,P)$ can also be
described as  the number of integral points in some polytope in an
Euclidean space. See Lemma \ref{come}.

\medskip
 As we will see explicitly in Section \ref{Sec Pieri rule for GL},
the combinatorial yoga that has been constructed to describe
representations of $\GL_n$ uses Young diagrams and semistandard
tableau. This formalism has also been extended to deal with the
classical isometry groups. This yoga can be used to describe the
Pieri rules
for isometry groups, but tensor products
for the isometry groups do not sit as comfortably in this formalism
as does branching. Our results show that  (at least under some
restrictions, which we refer to as the {\it stable range condition})
there is a uniform type of structure, that of Hibi rings, that can
describe tensor products as well as branching for all classical
groups. Moreover, the posets for the Hibi rings of iterated Pieri
rules have a family resemblance to the Hibi rings of Hodge's
standard monomial theory. This raises the possibility that there is
a uniform theory that encompasses all iterated branching algebras
and all iterated Pieri algebras for all classical groups. We hope to
study this issue in a future publication.

\medskip
This paper is arranged as follows: In  Section 2, we introduce
notation for the representations of $\Gln$, $\On$ and $\Sptn$, and
review the Pieri rule for $\Gln$. We define the Pieri algebra
$\faol$ in Section 3, where $2(k+\ell)<n$. In Section 4, we
compute the dimension of the homogeneous components of $\faol$.
  We study the
structure of $\faol$ in Section 5. Finally, we briefly explain the
parallel theory for $\Sptn$ in Section 6.

\bigskip\noindent {\bf Acknowledgement:} We thank Roger
Howe for helpful conversations.

\bigskip
\section{Preliminaries} In this section,
we introduce notation for the representations of $\Gln$, $\On$ and
$\Sptn$, and review the Pieri rule for $\Gln$. We will use the following notation throughout this paper: If $G$ is a group and $V$ a $G$ module, then $V^G$ will denote the space of all vectors in $V$ fixed by $G$.

\subsection{Representations of $\Gln$}
Let $\Gln=\Gln(\C)$ denote the general linear group consisting of
all $n\times n$ invertible complex matrices, and let $B_n=A_n U_n$
be the standard Borel subgroup of upper triangular matrices in
$\Glnn$, where $A_n$ is the diagonal torus in $\Glnn$ and $U_n$ is
the maximal unipotent subgroup consisting of all  the upper
triangular matrices with $1$'s on the diagonal.

\medskip  Recall that a {\em Young diagram} $D$ is
an array of square boxes arranged in left-justified horizontal rows,
with each row no longer than the one above it (\cite{Fu2}). If $D$
has at most $m$ rows, then we shall write it as
\[D=(\lambda_1,...,\lambda_m)\]
where for each $i$, $\lambda_i$ is the number of boxes in the $i$-th
row of $D$. We shall denote the number of rows in $D$ by $r(D)$, and
$|D|=\lambda_1+\cdots+\lambda_m$.

\medskip
For a
Young diagram $D=(\lambda_1,...,\lambda_n)$ with at most $n$ rows, let
$\psi^{D}_n:A_n\rightarrow\C^\times$ be the character given by
\begin{equation}\label{psid}
\psi^{D}_n[{\mathrm{diag}}(a_1,...,a_n)]=a^{\lambda_1}_1
\cdot\cdot\cdot a^{\lambda_n}_n.
\end{equation}
Let $\hat{A}^+_n$ be the semigroup of dominant weights for $\Glnn$
with respect to the Borel subgroup $B_n$. Then $\psi^{D}_n\in
\hat{A}^+_n$ (\cite{GW}), and  we shall denote the irreducible
representation of $\Gln$  with highest weight $\psi^{D}_n$ by
$\rho^D_n$.  By the theory of highest weight (\cite{GW}), the space
$\left(\rho^D_n\right)^{U_n}$ of $U_n$-invariants in $\rho^D_n$ is
one-dimensional, and $A_n$ acts on $\left(\rho^D_n\right)^{U_n}$ by the character $\psi^D_n$, i.e. the nonzero elements in
$\left(\rho^D_n\right)^{U_n}$ are the $\Gln$ highest weight vectors of
weight $\psi^D_n$.

\subsection{Representations of $\On$}\label{son}
Let $\On=\On(\C)$ be the subgroup of $\Glnn$ which preserves the
symmetric bilinear form

\begin{equation}\label{onform}
\langle \begin{pmatrix}u_1,\\ \vdots\\u_n\end{pmatrix},
\begin{pmatrix}v_1,\\ \vdots\\v_n\end{pmatrix}
 \rangle=\sum_{j=1}^nu_jv_{n-j+1}
\end{equation}
on $\C^n$. The irreducible finite dimensional representations of
$\On$ are parameterized by Young diagrams $D$ such that the sum of
the lengths of the first two columns of $D$ does not exceed $n$
(\cite{Wy, GW, Ho1}). For such a Young diagram $D$, we
shall denote the $\On$ representation associated with $D$ by
$\sigma^D_n$. Specifically, $\sigma^D_n$ is the irreducible
representation of $\On$ generated by the $\Gln$ highest weight
vector in $\rho^D_n$. See Section 3.6 of \cite{Ho1} for more
details.

\medskip
Let $\SO$ denote the subgroup of $\On$ containing elements of $\On$
with determinant $1$, and let
\[\AS=A_n\cap\SO,\hs \Nn=U_n\cap \SO.\]
Explicitly,
\[\AS=\left\{\begin{array}{ll}
\{\mathrm{diag}(a_1,...,a_m,a_m^{-1},...,a_1^{-1}): \
a_1,...,a_m\in\C^\times\}&n=2m\\
\{\mathrm{diag}(a_1,...,a_m,1,a_m^{-1},...,a_1^{-1}): \
a_1,...,a_m\in\C^\times\}&n=2m+1.\\
\end{array}\right.\]

\medskip If $D$ is a Young diagram with $r(D)\neq n/2$, then the restriction of $\sigma^D_n$ to $\SO$ is
irreducible.  If in addition, $r(D)<n/2$ and
$\phi^D_n:A_{\mathrm{SO}_n} \rightarrow \C^\times$ is the
restriction of the character $\psi^D_n$ to $\AS$, then as a $\SO$
module, $\sigma^D_n$ has highest weight $\phi^D_n$. In this case, we
also have $(\sigma^D_n)^{\Nn}=(\rho^D_n)^{U_n}$.

\subsection{Representations of $\Sptn$}\label{spr}
Let $\Sptn=\Sptn(\C)$ be the
 subgroup of
$\GL_{2n}$ which preserves  the symplectic form $(.,.)$ on $
\C^{2n}$  given by
\begin{equation}\label{spform}
( \begin{pmatrix}
x_1\\ \vdots \\ x_n\\y_1\\ \vdots\\y_n
 \end{pmatrix},
\begin{pmatrix}
x^\prime_1\\ \vdots\\ x^\prime_n\\y^\prime_1\\ \vdots\\y^\prime_n
\end{pmatrix}
 ) =\sum_{j=1}^n(x_jy^\prime_{j}-y_jx^\prime_j).
\end{equation}
The diagonal torus $A_{\Sptn }$ of $\Sptn$ is isomorphic to
$(\C^\times)^n$. The highest weights and hence the irreducible
finite dimensional representations of $\Sptn$ are parametrized by
Young diagrams with at most $n$ rows. If $D$ is such a Young
diagram, we shall denote the corresponding highest weight and
representation by $\chi^D_{2n}$ and $\tau^D_{2n}$ respectively. See
Section 3.8 of \cite{Ho1} for more details.

\medskip
\subsection{Pieri rule for $\Gln$ and Kostka coefficients}\label{Sec Pieri rule for GL}
The Pieri rule for $\Gln$ describes the decomposition of a tensor
product of a general irreducible representation with a
representation corresponding to a Young diagram with only one row.
In this section, we shall review this rule and one of its
generalization.

\medskip\noindent {\bf Definition.} For two Young diagrams $A=(a_{1},a_{2},\cdots )$ and
$B=(b_{1},b_{2},\cdots ) $, we denote $A\sqsupseteq B$ (or
$B\sqsubseteq A$)  and say $A$ \textit{interlaces} $B$ if
\[
a_{j}\geq b_{j}\geq a_{j+1}\ \ \ \ \mbox{for all $j$.}
\]

\medskip\nnn The following result is well known (\cite{GW, Ho1}):

\medskip \nnn {\bf Pieri rule for $\Gln$.} {\em Let $D$ be a Young
diagram with at most $n$ rows and let $p$ be a nonnegative integer.
Then
\[\rho^D_n\otimes \rho^{(p)}_n=\sum_{ F}\rho^F_n\]
where the sum is taken over all Young diagrams $F$ with at most $n$
rows such that $D\sqsubseteq F$ and $|F|-|D|=p$.}

\medskip By applying the Pieri rule repeatedly, we can give a
description of the decomposition of the tensor product
\begin{equation}\label{tensor}
\rho^D_n\otimes \rho^{(p_1)}_n\otimes\cdots\otimes\rho^{(p_\ell)}_n.
\end{equation}
In fact, a representation $\rho^F_n$ occurs in this tensor product
if and only if there is a sequence of Young diagrams
$(F_0,F_1,...,F_{\ell})$ such that \begin{equation}\label{chain}
D=F_0\sqsubseteq F_1\sqsubseteq\cdots\sqsubseteq
F_{\ell-1}\sqsubseteq F_\ell=F,
\end{equation}
and \begin{equation}\label{nb} |F_j|-|F_{j-1}|=p_j, \hs 1\leq j\leq
\ell.\end{equation}
Moreover, the number of sequences
$(F_0,F_1,...,F_{\ell})$ satisfying \eqref{chain} and \eqref{nb}
 gives the multiplicity of
$\rho^F_n$ in the tensor product in \eqref{tensor}.

\medskip
We shall give another description of this multiplicity. If a  Young
diagram $D$ sits inside another Young diagram $F$, then we write
$D\subset F$. In this case, by removing all boxes belonging to $D$,
we obtain the {\em skew diagram} $F/D$. If we put a positive number
in each box of $F/D$, then it becomes a {\em skew tableau} and we
say that the {\em shape} of this skew tableau is $F/D$. If the
entries of this skew tableau are taken from $\{1,2,...,m\}$, and
$\mu_j$ of them are $j$ for $1\leq j\leq m$, then we say the {\em
content} of this skew tableau is $E=(\mu_1,...,\mu_m)$. A skew
tableau $T$  is called {\em semistandard} if the numbers in each row
of $T$  weakly increase from left-to-right, and the numbers in each
column of $T$ strictly increase from top-to-bottom. The number of
semistandard tableaux of shape $F/D$ and with content $E$ is denoted
by $K_{F/D,E}$ and it is called a {\em skew Kostka number}
(\cite{Sta}).

\medskip We now fix two Young diagrams $D\subset F$, and consider a
sequence $(F_1,...,F_{\ell})$ which satisfies conditions
\eqref{chain} and \eqref{nb}.  We regard $F/D$  as a union of the
skew diagrams $F_i/F_{i-1}$, $1\leq i\leq\ell$. By filling the boxes
in $F_i/F_{i-1}$ with $i$ for each $1\leq i\leq \ell$, we obtain a
semistandard tableau $T$ of shape $F/D$ and content
$P=(p_1,...,p_\ell)$. This sets up a bijection between sequences of
Young diagrams $(F_1,...,F_{\ell})$ which satisfy conditions
\eqref{chain} and \eqref{nb}, and semistandard tableaux of shape
$F/D$ and content $P$. It follows that $K_{F/D,P}$ is the
multiplicity of $\rho^F_n$ in the tensor product \eqref{tensor}.
This proves:

\medskip \nnn {\bf Iterated Pieri rule for $\Gln$.}
{\em Let $D$ be a Young diagram with at most $n$ rows and
$P=(p_1,...,p_\ell)$ be a sequence of nonnegative integers. Then
\begin{equation}\label{tensor1}
\rho^D_n\otimes
\rho^{(p_1)}_n\otimes\cdots\otimes\rho^{(p_{\ell})}_n
=\sum_FK_{F/D,P}\rho^F_n.
\end{equation}}

\section{The construction of a Pieri algebra for $\On$}\label{sfaol}
Consider a tensor product of $\On$  representations of the form
$\sigma^D_n\otimes\sigma^E_n$ where the Young diagram $E$ consists
of only one row, that is, $E=(p)$ for some nonnegative integer $p$.
Under the diagonal action of $\On$,  $\sigma^D_n\otimes\sigma^E_n$
is decomposed as a direct sum
\begin{equation}\label{onef}
\sigma^D_n\otimes\sigma^E_n=\bigoplus_Fm_F\sigma^F_n
\end{equation}
where for each Young diagram $F$ appearing in the sum, $m_F$ is the
multiplicity of $\sigma^F_n$ in $\sigma^D_n\otimes\sigma^E_n$. A
description of the multiplicity $m_F$ is called the {\em Pieri Rule}
for $\On$.

\medskip It is natural to generalize  \eqref{onef} to the
tensor product of $\sigma^D_n$ with any number of representations
indexed by one-row Young diagrams:
\[\sigma^D_n\otimes\sigma^{(p_1)}_n
\otimes\sigma^{(p_2)}_n\otimes\cdots \otimes\sigma^{(p_\ell)}_n.\]
We call a description of the multiplicities in this tensor product
the {\em iterated Pieri rule} for $\On$. In this section, we shall
construct an algebra, called a {\em Pieri algebra} of $\On$, whose
algebra structure encodes information on the iterated Pieri rule for
the case when $r(D)\leq k$ and $2(k+\ell)<n$. We call the condition
$2(k+\ell)<n$ the {\em stable range condition}.

\subsection{A realization of irreducible representations of $\On$}
\label{relon} Let $n$ and $k$ be positive integers such that $2k<n$,
and let $\Mnk=\Mnk(\C)$ denote the space of all $n\times k$ complex
matrices. Let $\pmnk$ be the algebra of polynomial functions on
$\Mnk$, that is, each $p\in\pmnk$ is of the form
\[p(x)=\sum_{\alpha}a_\alpha x^\alpha,\] where
the sum is finite, $x=(x_{ij})\in\Mnk$, each $\alpha=(\alpha_{ij})$
appearing in the sum is an $n\times k$ matrix of nonnegative
integers, $a_\alpha\in\C$ and
\[x^\alpha=\prod_{i,j}x^{\alpha_{ij}}_{ij}.\]
We define an action of $\On\times\Glk$ on $\pmnk$  as follows: For
$(g,h)\in \On\times\Glk$, $f\in \pmnk$ and $X\in \Mnk$, let
\begin{equation}\label{Onaction}
\left((g,h).f\right)(X)=f(g^{-1} X h).
\end{equation}
 For $1\leq i,j\leq k$, let
\[r_{i,j}(X)=\langle X_i,X_j\rangle,\]
where $\langle.,.\rangle$ is the symmetric bilinear form defined in
\eqref{onform} and  $X_i$ and $X_j$ are the $i$-th and the $j$-th columns of $X$
respectively. These polynomials generate the algebra of $\On$
invariants in $\pmnk$. Let $I_{nk}$ be the ideal of $\pmnk$
generated by $\{r_{i,j}:\ 1\leq i,j\leq k\}$. It is stable under the
action by $\On\times \Glk$. So the action \eqref{Onaction} induces
an action by $\On\times\Glk$ on  the quotient algebra
$\pmnk/I_{nk}$. Under this action, $\pmnk/I_{nk}$  has a
decomposition given by (\cite{Ho1})
\begin{equation}\label{quotientdecom}
\pmnk/I_{nk} \cong\sum_{r(D)\leq k}\sigma^D_n\otimes\rho^D_k.
\end{equation}
Let $\left( \pmnk/I_{nk} \right)^{U_k}$ be the subalgebra of
$\pmnk/I_{nk}$ consisting of all the elements fixed by the maximal
unipotent subgroup $U_k$ of $\Glk$. Then by extracting $U_k$
invariants in \eqref{quotientdecom}, we obtain
\[\left( \pmnk/I_{nk} \right)^{U_k} \cong\sum_{r(D)\leq
k}\sigma^D_n\otimes \left(\rho^D_k\right)^{U_k}.\]   By the theory of highest
weight, $\dim \left(\rho^D_k\right)^{U_k}=1$ and the diagonal torus
$A_k$ of $\Glk$ acts on $\left(\rho^D_k\right)^{U_k}$ by the
character $\psi^D_k$   given in \eqref{psid}. It follows that $\left( \pmnk/I_{nk}
\right)^{U_k}$ is a module for $\On\times A_k$. For each Young
diagram $D$ with $r(D)\leq k$,
\[\sigma^D_n\cong \sigma^D_n\otimes \left(\rho^D_k\right)^{U_k}.\]
Thus we can realize the irreducible representation $\sigma^D_n$   as the $\psi^D_k$-eigenspace of $A_k$ in $\left( \pmnk/I_{nk}
\right)^{U_k}$.

\medskip We also need to realize those $\On$ representations indexed by
one-row Young diagrams. To do this, we simply repeat the above
construction with $k=1$.  For each $1\leq j\leq\ell$,  let $\C^n_j$
be a copy of $\C^n$. We shall denote a typical vector in
$\C^n_j$ as \[Y_j=\begin{pmatrix}y_{1,j}\\y_{2,j}\\
\vdots\\y_{n,j}\end{pmatrix},\] so the algebra $\calp(\C^n_j)$
of polynomial functions on $\C^n_j$ can be regarded as the
polynomial algebras on the variables $y_{1,j},y_{2,j},...,y_{n,j}$.
Let
\[t_{j}(Y_j)=\langle Y_j,Y_j\rangle,\] and let $I^{(j)}_{n1}$ be the
ideal of $\calp(\C^n_j)$ generated by $t_j$. We form the quotient
algebra $\calp(\C^n_j)/I^{(j)}_{n1}$. Then
$\calp(\C^n_j)/I^{(j)}_{n1}$ is a module for $\On\times \GL_1$ and
it can be decomposed as
\[\calp(\C^n_j)/I^{(j)}_{n1}\cong
\sum_{p_j\geq 0}\sigma^{(p_j)}_n\otimes\rho^{(p_j)}_1.\] As before,  the
$\On$ representation $\sigma^{(p_j)}_n$ can now be identified with
the $\psi^{(p_j)}_1$-eigenspace of $\GL_1$ in
$\calp(\C^n_j)/I^{(j)}_{n1}$.

\subsection{The algebra $\faol$}
We shall assume $2(k+\ell)<n$ from now on. To realize the tensor
product $\sigma^D_n\otimes \sigma^{(p_1)}_n\otimes\cdots\otimes
\sigma^{(p_\ell)}_n$ of $\On$ representations, we form the tensor
product of the algebras
\[\fto:=\left( \pmnk/I_{nk} \right)^{U_k}\otimes
 \left(\calp(\C^n_1)/I^{(1)}_{n1}\right)\otimes\cdots\otimes
 \left(\calp(\C^n_\ell)/I^{(\ell)}_{n1}\right).\]
It is a module for $\On\times A_k\times \GL^{\ell}_1$, and it can be
decomposed as
\begin{eqnarray*}  \fto
  &\cong& \left\{\sum_{r(D)\leq k}\sigma^D_n\otimes
\left(\rho^D_k\right)^{U_k}\right\}\otimes \left\{\sum_{p_1\geq
0}\sigma^{(p_1)}_n\otimes \rho^{(p_1)}_1\right\}\otimes
\cdots\otimes\left\{\sum_{p_\ell\geq 0}\sigma^{(p_\ell)}_n\otimes
\rho^{(p_\ell)}_1\right\}\\
&\cong&  \sum_{\stackrel{\scriptstyle r(D)\leq k}{p_1,...,p_\ell\geq
0}}\left(\sigma^D_n\otimes \sigma^{(p_1)}_n\otimes\cdots\otimes
\sigma^{(p_\ell)}_n\right)\otimes \left(
\rho^D_k\right)^{U_k}\otimes \rho^{(p_1)}_1\otimes\cdots\otimes
\rho^{(p_\ell)}_1.
\end{eqnarray*}
We can identify $\sigma^D_n\otimes
\sigma^{(p_1)}_n\otimes\cdots\otimes \sigma^{(p_\ell)}_n$ with the
$\psi^D_k\times\psi^{(p_1)}_1\times\cdots\times\psi^{(p_\ell)}_1$-eigenspace
of $A_k\times \GL^\ell_1$ in $\fto$. To decompose this tensor
product, we extract the $\Nn$ invariants from this space. Let
\begin{equation}\label{dfaol}
\faol:=\left(\fto\right)^{\Nn}, \end{equation} that is, $\faol$ be
the subalgebra of $\fto$ consisting of all elements of $\fto$ fixed
by $\Nn$. It is a module for $\AS\times A_k\times \GL_1^\ell$ and
can be decomposed as
\[ \faol\cong
\sum_{\stackrel{\scriptstyle r(D)\leq k}{p_1,...,p_\ell\geq
0}}\left(\sigma^D_n\otimes \sigma^{(p_1)}_n\otimes\cdots\otimes
\sigma^{(p_\ell)}_n\right)^\Nn\otimes \left(
\rho^D_k\right)^{U_k}\otimes \rho^{(p_1)}_1\otimes\cdots\otimes
\rho^{(p_\ell)}_1.\] The space $\left(\sigma^D_n\otimes
\sigma^{(p_1)}_n\otimes\cdots\otimes \sigma^{(p_\ell)}_n\right)^\Nn$
can be further decomposed as a sum of eigenspaces of $\AS$, so we
can write \begin{equation}\label{homcom}
\faol=\sum_{F,D,P}
 \cale_{F,D,P},
\end{equation}
 where the sum is taken over all Young diagrams $D$ and $F$ with
$r(D)\leq k$, $r(F)\leq k+\ell$ and all
$P=(p_1,...,p_\ell)\in\Z^\ell_{\geq 0}$, and $\cale_{F,D,P}$ is the
$\phi^F_n\times\psi^D_k\times\psi^{(p_1)}_1\times\cdots\times
\psi^{(p_\ell)}_1$-eigenspace of $\AS\times A_k\times \GL_1^\ell$
(recall that $\phi^F_n$ is the restriction of $\psi^F_n$ to $\AS$).
Since $\AS\times A_k\times \GL^\ell_1$ acts on $\faol$ by algebra
automorphisms, the direct sum decomposition \eqref{homcom} defines a
multi-grading on $\faol$.

\medskip Observe that the homogeneous component $\cale_{F,D,P}$ can be identified
with the space of all $\SO$ highest weight vectors of weight
$\phi^F_n$ in $\sigma^D_n\otimes
\sigma^{(p_1)}_n\otimes\cdots\otimes \sigma^{(p_\ell)}_n$.
Now the Young diagrams which label the irreducible $\On$ representations occurring  in the
tensor product $\sigma^D_n\otimes
\sigma^{(p_1)}_n\otimes\cdots\otimes \sigma^{(p_\ell)}_n$ 
have at most $k+\ell$ rows. Since $k+\ell<n/2$, these $\On$ representations
 remain irreducible under the action by $\SO$ and they are
determined by the $\SO$ highest weight vectors they contain.
Consequently, the dimension of $\cale_{F,D,P}$ coincides with the
multiplicity of $\sigma^F_n$ in $\sigma^D_n\otimes
\sigma^{(p_1)}_n\otimes\cdots\otimes \sigma^{(p_\ell)}_n$. Thus the
algebra structure of $\faol$ encodes the iterated Pieri rule in the
case $2(k+\ell)<n$. In view of this property, we will call $\faol$
an {\em   $\On$ Pieri algebra.} The remaining paper will be devoted
to determining the structure of this algebra.

\medskip\section{The generalized Pieri rule for $\On$}\label{secfour}
In this section, we show that the $\On$ Pieri algebra $\faol$ is
isomorphic to a subalgebra   of a polynomial algebra. Using this
isomorphism, we determine the dimension of the homogeneous
components $\cale_{F,D,P}$ of $\faol$.

\medskip As before, we assume that $2(k+\ell)<n$.
The group $\On\times\GL_{k+\ell}$ acts on the algebra
$\calp(\rM_{n,k+\ell})$ of polynomial functions on $\rM_{n,k+\ell}$
in a similar way as in Section \ref{relon} (with $k$ replaced by
$k+\ell$), and we restrict this action to the subgroup $\On\times
\left(\Glk\times\GL_1^\ell\right)$ of $\On\times\GL_{k+\ell}$ . We
write a typical matrix in $\rM_{n,k+\ell}$ as $X=(x_{ij})$, and the
$j$th column of $X$ as $X_j$. Let
\[r_{ij}=\langle X_i,X_j\rangle\hs (1\leq i,j\leq k+\ell),\]
 \[R=\{r_{ij}:\ 1\leq
i,j\leq k\}\cup\{r^2_{k+j,k+j}:\ 1\leq j\leq \ell\},\] and let $J$
be the ideal of $\pmnkl$ generated by $R$. Then $J$ is stable under
under the action by $\On\times \left(\Glk\times\GL_1^\ell\right)$,
so the action by $\On\times \left(\Glk\times\GL_1^\ell\right)$ on
$\pmnkl$ induces an action on the quotient algebra $\pmnkl/J$. There
is an obvious   isomorphism of algebras and $\On\times
\left(\Glk\times\GL_1^\ell\right)$ modules:
\begin{equation}\label{tensoriso}
\left( \pmnk/I_{nk} \right)\otimes
\left(\calp(\C^n_1)/I^{(1)}_{n1}\right)%
\otimes\cdots\otimes
\left(\calp(\C^n_\ell)/I^{(\ell)}_{n1}\right)\nonumber
\newline
\cong \calp(\rM_{n,k+\ell})/J.
\end{equation}
Here, we identify $X_{k+j}=Y_j$, i.e. $x_{i,k+j}=y_{i,j}$, $1\leq
i\leq n$ and $1\leq j\leq \ell$, so that $r_{k+j,k+j}=t_j$ (see
Section \ref{sfaol} for notation). It now follows from this and the
definition of $\faol$ given in equation \eqref{dfaol} that
\begin{equation}\label{ptensorisopmnkl}
\faol\cong \left(\calp(\rM_{n,k+\ell})/J\right)^{\Nn\times U_k}
\end{equation}
as $\AS\times A_k\times \GL^\ell_1$ modules and algebras.

\medskip\begin{prop} \label{ispa} The $\On$ Pieri algebra $\faol$ is
isomorphic to a subalgebra of a polynomial algebra.
\end{prop}

\proof For $1\leq i,j\leq k+\ell$, let
 \[ \Delta_{ij}=\sum_{a=1}^n\frac{\partial^2}{\partial
x_{a,i}\partial
 x_{n-a+1,j}},\]
and let
\[\calh(\rM_{n,k+\ell})=\{f\in \calp(\rM_{n,k+\ell}):\
\Delta_{ij}(f)=0\ \forall 1\leq i\leq j\leq k+\ell\}\] be the space
of $\On$ harmonic polynomials in $\pmnkl$.  Let $\pmnkl^{\On}$ be
the subalgebra of $\On$ invariant polynomials on $\rM_{n,k+\ell}$,
which is a polynomial algebra on $r_{ij}$, $1\leq i\leq j\leq k$.
Both $\calh(\rM_{n,k+\ell})$ and $\pmnkl^{\On}$ are stable under the
action by $\On\times\GL_{k+\ell}$, so the tensor product
$\calhmkl\otimes \pmnkl^{\On}$ is an $\On\times\GL_{k+\ell}$ module.
We consider the multiplication map
\[m:\calhmkl\otimes \pmnkl^{\On}\rightarrow \pmnkl,\]
\begin{equation}\label{hq}
m(h\otimes q)=hq.\end{equation} Since $2(k+\ell)<n$, $m$ is an
$\On\times\GL_{k+\ell}$ module isomorphism  (\cite{Ho1}).

\medskip We now label the various copies of
$\GL_1$ in $\Glk\times \GL_1^\ell$ as $\GL^{(j)}_1$,
$j=1,2,...,\ell$, and note that for each $1\leq j\leq \ell$, the
subalgebra of $\pmnkl^{\On}$ generated by $\{r_{1,k+j}, ...,
r_{k,k+j}\}$ is stable under the action by $(\Glk\times\GL^{(j)}_1)$
and is isomorphic to $\calp(\C^k)$. We shall denote this subalgebra
as $\calp(\C^k_j)$, where $\C^k_j$ denotes a copy of $\C^k$.
Similarly, for  $ 1\leq s<t\leq \ell$, the subalgebra of
$\pmnkl^{\On}$ generated by $r_{k+s,k+t}$ is stable under the action
by $\GL^{(s)}_1\times\GL^{(t)}_1$, and is isomorphic to $\calp(\C)$.
So we shall denote this subalgebra by $\calp(\C_{s,t})$, where
$\C_{s,t}$ denote a copy of $\C$.  The map $m$ now induces
an isomorphism of $\On\times(\Glk\times \GL_1^\ell)$ modules:
\begin{equation}
\pmnkl/J\cong \calhmkl\otimes \left(\bigotimes_{j=1}^\ell \calp(\C%
^k_j)\right)\otimes\left(\bigotimes_{1\leq s<t\leq \ell} \calp(\C%
_{s,t})\right).  \label{decom1}
\end{equation}
Using this and  \eqref{ptensorisopmnkl},  we obtain an
$\AS\times(A_k\times \GL_1^\ell)$ module isomorphisms:
 \begin{eqnarray*}
\faol
&\cong&\left\{ \calhmkl\otimes \left(\bigotimes_{j=1}^\ell \calp(\C%
^k_j)\right)\otimes\left(\bigotimes_{1\leq s<t\leq \ell} \calp(\C%
_{s,t})\right) \right\}^{\Nn\times U_k}    \\
&=& \left\{ \calhmkl^\Nn\otimes \left(\bigotimes_{j=1}^\ell \calp(\C%
^k_j)\right)\otimes\left(\bigotimes_{1\leq s<t\leq \ell} \calp(\C%
_{s,t})\right) \right\}^{U_k}   \\
&=& \left\{ \pmnkl^{U_n}\otimes \left(\bigotimes_{j=1}^\ell \calp(\C%
^k_j)\right)\otimes\left(\bigotimes_{1\leq s<t\leq \ell} \calp(\C%
_{s,t})\right) \right\}^{U_k}\\
&=&\left(\calp_{n,k,\ell}\right)^{U_n\times U_k},
\end{eqnarray*}
where $\calp_{n,k,\ell}$ is the polynomial algebra
\[ \pmnk\otimes \left(\bigotimes_{i=1}^\ell \calp(\C
^n_i)\right)\otimes \left(\bigotimes_{j=1}^\ell \calp(\C^k_j)\right)\otimes%
\left(\bigotimes_{1\leq s<t\leq \ell} \calp(\C_{s,t})\right),\]
and for each $1\leq i\leq \ell$,  $\C^n_i$ is a copy of $\C^n$.
In establishing the isomorphisms above, we have used the fact that $\calhmkl^\Nn=\pmnkl^{U_n}$.
Since this module isomorphism is induced by the multiplication map \eqref{hq}, the isomorphism for $\faol\cong \left(\calp_{n,k,\ell}\right)^{U_n\times U_k}$ is actually an algebra isomorphism.  $\Box$

\medskip In view of proposition \ref{ispa},
we shall identify $\faol$ with
$\left(\calp_{n,k,\ell}\right)^{U_n\times U_k}$ from now on.

\medskip\nnn{\bf Notation.}
For $A=(a_i),B=(b_j)\in\R^\ell$ and $C=(c_{st})_{1\leq s<t\leq
\ell}\in \R^{\ell(\ell-1)/2}$, let $S(A,B,C)=(p_i)\in\R^\ell$ be
defined by
\begin{equation}\label{pcon}
 p_i=a_i+b_i+c_{1i}+\cdots+c_{i-1,i}+c_{i,i+1}+\cdots+c_{i,\ell}\ \
 (1\leq i\leq \ell).
 \end{equation}

\medskip
\begin{prop}\label{dimf} The dimension of the homogeneous component
$\cale_{F,D,P}$ is given by
\[\dim \cale_{F,D,P}=
 \sum_{\stackrel{r(E)\leq
k}{S(A,B,C)=P}}K_{F/E,A} K_{D/E,B}.
\]
\end{prop}\proof   By
$(\Gln,\GL_{m})$ duality (\cite{Ho1}), we have
\begin{equation*}
\pmnk\cong \sum_{r(E)\leq k}\rho _{n}^{E}\otimes \rho _{k}^{E},\ \ \ \ \calp(%
\C_{s,t})\cong \sum_{c_{s,t}\geq 0}\rho _{1,s}^{(c_{s,t})}\otimes
\rho _{1,t}^{(c_{s,t})}\ \ (1\leq s<t\leq \ell ),
\end{equation*}%
\begin{equation*}
\calp(\C_{i}^{n})\cong \sum_{a_{i}\geq 0}\rho _{n}^{(a_{i})}\otimes
\rho _{1,i}^{(a_{i})},\ \ \ \ \calp(\C_{j}^{k})\cong \sum_{b_{j}\geq
0}\rho _{k}^{(b_{j})}\otimes \rho _{1,j}^{(b_{j})}\ \ (1\leq i,j\leq
\ell ).
\end{equation*}
Here, to indicate that a particular copy $\GL^{(j)}_1$ of $\GL_1$ acts on the
representation $\rho^{(m)}_1$, we have written $\rho^{(m)}_1$ as
$\rho^{(m)}_{1,j}$. Then \begin{eqnarray*}&&\calp_{n,k,\ell}=
\pmnk\otimes \left(\bigotimes_{i=1}^\ell \calp(\C^n_i)\right)\otimes \left(\bigotimes_{j=1}^\ell \calp(\C^k_j)\right)\otimes
\left(\bigotimes_{1\leq s<t\leq \ell} \calp(\C_{s,t}) \right)\\
 &\cong &\sum_{\overset{\scriptstyle r(E)\leq
k}{A,B,C}}\left\{ \left( \rho _{n}^{E}\otimes \left(
\bigotimes_{i=1}^{\ell }\rho _{n}^{(a_{i})}\right) \right)
 \otimes \left( \rho _{k}^{E}\otimes \left(
\bigotimes_{j=1}^{\ell }\rho _{k}^{(b_{j})}\right) \right)
 \right.  \\
&&\hspace{1in}\left. \otimes \left( \bigotimes_{1\leq s<t\leq \ell
}\rho _{1,s}^{(c_{s,t})}\otimes \rho _{1,t}^{(c_{s,t})}\right)
\otimes \left( \bigotimes_{i=1}^{\ell }\rho _{1,i}^{(a_{i})}\right)
\otimes \left(
\bigotimes_{j=1}^{\ell }\rho _{1,j}^{(b_{j})}\right) \right\}  \\
&\cong &\sum_{\overset{\scriptstyle r(E)\leq k}{A,B,C}}\left\{
\left(\sum_F K_{F/E,A}\ \rho _{n}^{F}\right) \otimes \left( \sum_D
K_{D/E,B}\ \rho _{k}^{D}\right) \otimes \left( \bigotimes_{1\leq
i\leq \ell }\rho
_{1,i}^{(p_{i})}\right) \right\}  \\
&\cong &\sum_{F,D,P}\left( \sum_{\overset{\scriptstyle r(E)\leq k}{
S(A,B,C)=P}} K_{F/E,A}K_{D/E,B}
 \right) \rho _{n}^{F}\otimes \rho
_{k}^{D}\otimes \left( \bigotimes_{1\leq i\leq \ell }\rho
_{1,i}^{(p_{i})}\right)
\end{eqnarray*}%
where $A=(a_1,...,a_\ell)$, $B=(b_1,...,b_\ell)$, $C=(c_{s,t})$ and
$P=S(A,B,C)=(p_1,...,p_\ell)$ in the sums. The proposition now
follows from extracting the $U_n\times U_k$ invariants from
$\calp_{n,k,\ell}$.
 $\Box $

\medskip
\begin{cor}
Let $k,\ell$ and $n$ be positive integers such that $2(k+\ell)<n$,
$D$ a Young diagram such that $r(D)\leq k$ and
$P=(p_1,...,p_\ell)\in\Z^\ell_{\geq 0}$. Then for any Young diagram
$F$ which labels a representation of $\On$, the multiplicity
$m(F,D,P)$ of $\sigma^F_n$ in the tensor product
\[\sigma^D_n\otimes\sigma^{(p_1)}_n
\otimes\sigma^{(p_2)}_n\otimes\cdots \otimes\sigma^{(p_\ell)}_n\]
    is given by
\[
m(F,D,P) =  \sum_{\stackrel{\scriptstyle r(E)\leq
k}{S(A,B,C)=P}}K_{F/E,A} K_{D/E,B}.\]
\end{cor}

\section{The structure of $\faol$}\label{stfaol}
In this section, we shall determine the structure of the algebra $\faol$. We first define a poset $\tG (k,\ell)$ which arises
from a consideration of the dimension formula given in Proposition
\ref{dimf}. The poset $\tG (k,\ell)$ gives rise to a finite
distributive lattice $\calg(k,\ell)$ as well as a generating set
$S(k,\ell)$ for the algebra $\faol$. By defining an appropriate
partial ordering on  $S(k,\ell)$, standard monomials on $S(k,\ell)$
form a basis for $\faol$. This standard monomial basis allows us to
show that $\faol$ has a flat deformation to the Hibi algebra
associated with the distributive lattice $\calg(k,\ell)$
(\cite{Hi, Ho2}).

\subsection{The poset $\Gamma(k,\ell)$} Suppose that
  the homogeneous component $\cale_{F,D,P}$ of $\faol$ is nonzero. Then by Proposition \ref{dimf},
\[\dim \cale_{F,D,P}=
 \sum_{\stackrel{r(E)\leq
k}{S(A,B,C)=P}}K_{F/E,A} K_{D/E,B}.
\]
Fix  $A=(a_1,...,a_\ell),B=(b_1,...,b_\ell)\in \Z^\ell_{\geq 0}$ and
$C=(c_{st})_{1\leq s<t<\ell}\in\Z^{\ell(\ell-1)/2}_{\geq 0}$ in the
sum. The Kostka number $K_{F/E,A}$ counts the number of sequences of
Young diagrams $(F_0,F_1,...,F_\ell)$ such that
\begin{equation}\label{seqn1}
E=F_0\sqsubseteq F_1\sqsubseteq\cdots\sqsubseteq
F_{\ell-1}\sqsubseteq F_\ell=F,\ \ \mbox{and}\ \
|F_j|-|F_{j-1}|=a_j\ \ (1\leq j\leq \ell). \end{equation} This sequence of Young diagrams can be identified with the ``truncated"
Gelfand-Tsetlin pattern  {\scriptsize
\[\begin{array}{ccccccccccccccccc }
\alpha_{\ell,1} && \alpha_{\ell,2} &&\alpha_{\ell,3}&&\cdot&&\cdot
&&\cdot&&\cdot&&\cdot  && \alpha_{\ell,k+\ell} \\
&\alpha_{\ell-1,1} && \alpha_{\ell-1,2} &&\alpha_{\ell-1,3}&& \cdot
&&\cdot&&\cdot&&\cdot && \alpha_{\ell-1,k+\ell-1}& \\
&&\cdot&&\cdot&&\cdot&&\cdot&&\cdot&&\cdot&&\cdot&&\\
&&&\cdot&&\cdot&&\cdot&&\cdot&&\cdot&&\cdot&& &\\
&&&&\alpha_{0,1}&&\alpha_{0,2}&&\cdot&&\cdot&&\alpha_{0,k}&&&&
\end{array}\]}

\noindent where $F_j=(\alpha_{j,1},\alpha_{j,2}, \cdots ,
\alpha_{j,k+j})$ for $0\leq j \leq \ell$, which is obtained from a Gelfand-Tsetlin
pattern with $\ell+k$ rows by removing its lowest $k-1$ rows.
Similarly, $K_{D/E,B}$ counts the number of sequences of Young
diagrams $(D_0,D_1,...,D_\ell)$ such that
\begin{equation}\label{seqn2}
E=D_0\sqsubseteq D_1\sqsubseteq\cdots\sqsubseteq
D_{\ell-1}\sqsubseteq D_\ell=D,\ \ \mbox{and}\ \
|D_j|-|D_{j-1}|=b_j\ \ (1\leq j\leq \ell), \end{equation} which can
be viewed as another ``truncated" Gelfand-Tsetlin pattern:
{\scriptsize
\[\begin{array}{ccccccccccccccc }
\beta_{\ell,1} &  & \beta_{\ell,2} &  &\beta_{\ell,3}&&\cdot&&\cdot   &  & \beta_{\ell,k} &  &  &  &  \\
& \beta_{\ell-1,1} &  & \beta_{\ell-1,2} &  &
\beta_{\ell-1,3}&&\cdot&&\cdot   &  & \beta_{\ell-1,k} &  &  &    \\
&&\cdot&&\cdot&&\cdot&&\cdot&&\cdot&&\cdot&&\\
&&&\cdot&&\cdot&&\cdot&&\cdot&&\cdot&&\cdot&\\
&  &  &  & \beta_{0,1} &  & \beta_{0,2} &  &
\beta_{0,3}&&\cdot&&\cdot   &  & \beta_{0,k}
\end{array}\]}

\noindent where $D_j=(\beta_{j,1},\beta_{j,2}, \cdots ,
\beta_{j,k})$ for $0\leq j \leq \ell$. Note that these two patterns
share the same lowest row, so we can invert the first pattern and
stack it below the second to obtain a pattern of the form
\begin{equation}\label{pattern}
\mbox{\scriptsize $\begin{array}{ccccccccccccccccc }
\beta_{\ell,1}&&\beta_{\ell,2}&&\cdot&&\cdot&&\beta_{\ell,k}&&&&&&&&\\
&\cdot&&\cdot&&\cdot&&\cdot&&\cdot&&&&&&&\\
&&\cdot&&\cdot&&\cdot&&\cdot&&\cdot&&&&&&\\
&&&\beta_{1,1}&&\beta_{1,2}&&\cdot&&\cdot&&\beta_{1,k}&&&&&\\
&&&&\alpha_{0,1}&&\alpha_{0,2}&&\cdot&&\cdot&&\alpha_{0,k}&&&&\\
&&&\cdot&&\cdot&&\cdot&&\cdot&&\cdot&&\cdot&& &\\
&&\cdot&&\cdot&&\cdot&&\cdot&&\cdot&&\cdot&&\cdot&&\\
&\alpha_{\ell-1,1} && \alpha_{\ell-1,2} &&\alpha_{\ell-1,3}&& \cdot
&&\cdot&&\cdot&&\cdot && \alpha_{\ell-1,k+\ell-1}& \\
\alpha_{\ell,1} && \alpha_{\ell,2} &&\alpha_{\ell,3}&&\cdot&&\cdot
&&\cdot&&\cdot&&\cdot  && \alpha_{\ell,k+\ell} \\
\end{array}$}\end{equation}
Thus the product $K_{F/E,A} K_{D/E,B}$ counts certain patterns of
this form. Following \cite{Ki1}, we can view the pattern
\eqref{pattern} as an order preserving function $f$ on a poset
$\Gamma(k,\ell)$, which is defined in the following way. Each entry
in \eqref{pattern} should correspond to a unique element in
$\Gamma(k,\ell)$, so that it is the image of this element  under the
function $f$. Moreover, the partial ordering on $\Gamma(k,\ell)$ is
induced by the partial ordering on the entries of the pattern. Thus
we let
\[\Gamma(k,\ell)=\{\gamma^{(j)}_i: \
1\leq i\leq \max(k,k+j),\ -\ell\leq j\leq \ell\},\] and the partial
ordering on $\Gamma(k,\ell)$ is given as follows: for $1\leq i\leq
\ell$,
\begin{equation*}
\gamma^{(-i)}_a \geq \gamma^{(-i+1)}_a \geq \gamma^{(-i)}_{a+1}\
\mbox{for}\ 1\leq a\leq k-1,\ \ \gamma^{(-i)}_k \geq
\gamma^{(-i+1)}_k,\ \mbox{and}
\end{equation*}
\begin{equation*}
\gamma^{(i)}_b \geq \gamma^{(i-1)}_b \geq \gamma^{(i)}_{b+1}\
\mbox{for}\ 1\leq b\leq k+i-1.
\end{equation*}
For example, if we display the elements of $\Gamma(k,\ell)$ as the
entries in \eqref{pattern}, then $\Gamma(k,1)$ and $\Gamma(k,2)$ are
given by
\begin{equation}
 \Gamma(k,1)=   \mbox{\scriptsize $
\begin{array}{cccccccccccccccc}
\gamma^{(-1)}_1 &  & \gamma^{(-1)}_2 &  & \gamma^{(-1)}_3 &  & \cdot
&  &
\cdot &  & \cdot &  & \gamma^{(-1)}_k &  &  &  \\
& \gamma^{(0)}_1 &  & \gamma^{(0)}_2 &  & \gamma^{(0)}_3 &  & \cdot
&  &
\cdot &  & \cdot &  & \gamma^{(0)}_k &  &  \\
\gamma^{(1)}_1 &  & \gamma^{(1)}_2 &  & \gamma^{(1)}_3 &  & \cdot &
& \cdot
&  & \cdot &  & \gamma^{(1)}_k &  & \gamma^{(1)}_{k+1} &  \\
&  &  &  &  &  &  &  &  &  &  &  &  &  &  &
\end{array}$}\label{Gamma ell=1}
\end{equation}

\bigskip
\[
\Gamma(k,2)=\mbox{\scriptsize $
\begin{array}{ccccccccccccc}
\gamma^{(-2)}_1 &  & \gamma^{(-2)}_2 &  & \cdot &  & \cdot &  &
\gamma^{(-2)}_k &  &  &  &  \\
& \gamma^{(-1)}_1 &  & \gamma^{(-1)}_2 &  & \cdot &  & \cdot &  &
\gamma^{(-1)}_k &  &  &  \\
&  & \gamma^{(0)}_1 &  & \gamma^{(0)}_2 &  & \cdot &  & \cdot &  &
\gamma^{(0)}_k &  &  \\
& \gamma^{(1)}_1 &  & \gamma^{(1)}_2 &  & \cdot &  & \cdot &  &
\gamma^{(1)}_k &  & \gamma^{(1)}_{k+1} &  \\
\gamma^{(2)}_1 &  & \gamma^{(2)}_2 &  & \cdot &  & \cdot &  &
\gamma^{(2)}_k
&  & \gamma^{(2)}_{k+1} &  & \gamma^{(2)}_{k+2} \\
&  &  &  &  &  &  &  &  &  &  &  &
\end{array}.$}
\]
The poset $\Gamma(k,\ell)$ is still insufficient for describing
$\dim\cale_{F,D,P}$; we also need to encode the information in the
matrix $C=(c_{s,t})_{1\leq s<t\leq \ell}$. We can view $C$ as a
function on another poset with $\ell(\ell-1)/2$ elements. Recall
that for posets $P$ and $Q$ on disjoint sets, the disjoint union
$P+Q$ is the poset on $P\cup Q$ such that $x\leq y$ in $P+Q$ if
either $x,y\in P$ and $x\leq y$ in $P$, or $x,y\in Q$ and $x\leq
y$ in $Q$. Let $F_\ell$ be the disjoint union of the singleton posets
$\{\ve_{s,t}\}$, $1\leq s<t\leq\ell$. Finally, we let
\begin{equation}\label{fposet}
\tG(k,\ell)=\Gamma(k,\ell)+F_\ell.
\end{equation}
The poset $\tG(k,\ell)$ will play a pivotal role in the structure of
$\faol$.

\medskip
\subsection{The lattice cone $\Omega(k,\ell)$} \label{okl}
By a {\em lattice cone}, we mean the intersection of a convex
polyhedral cone in $\R^n$ for some $n$ with $\Z^n$. It has a
structure of an {\em affine semigroup} with respect to vector
addition,
 i.e. it is a
finitely generated subsemigroup of $\Z^n$ containing $0$. In
\cite{Ho2}, Howe describes how to construct a lattice cone from a
finite poset, and gives the generators and relations for this
lattice cone. In this subsection, we shall use results in \cite{Ho2}
to describe the lattice cone $\Omega(k,\ell)$ corresponding to
$\tG(k,\ell)$.

\medskip We call a real-valued function $f$ on $\tG=\tG(k,\ell)$ \emph{%
order preserving} if
\begin{equation*}
x,y\in\tG,\ x \geq y \ \Longrightarrow \ f(x)\geq f(y).
\end{equation*}
Let
\[\calc=\calc(k,\ell)=\{f:\tG\rightarrow \R_{\geq 0}|
\ \mbox{$f$ is order preserving}\},\] and let
$\R^{\tG}=\R^{\tG(k,\ell)}$ be the space of all functions $ f:\tG
\rightarrow\R$. By listing the elements of $ \tG$ in a specific
order, say $u_1,u_2,...,u_N$ where
\begin{equation}\label{von}
N=(2\ell+1)k+\ell^2,
\end{equation}
we can identify each function $f\in \R^{\tG }$
  with the point
\begin{equation*}
(f(u_1),f(u_2),...,f(u_N))\in\R^N.
\end{equation*}
With this identification, it is easy to see that the subset of
$\R^N$ corresponding to $\calc $ is a convex polyhedral cone. Let
\begin{equation*}
\Omega=\Omega(k,\ell)=\calc\cap\Z^{\tG},
\end{equation*}
where $\Z^{\tG }$ is the space of all integer-valued functions on $
\tG$ (which can be identified with the lattice $\Z^N$ in $\R^N$). So
$\Omega $ is a lattice cone. In particular, it forms an affine
semigroup with the usual addition of functions.


\medskip Each element of
$\Omega$ gives rise to a pattern of the form \eqref{pattern} and a
matrix $C\in\Z^{\ell(\ell-1)}_{\geq 0}$.   To describe this more
precisely, we need some notation. If $f\in \calc$, we write
\begin{eqnarray} f_{i,j}&=&f\left( \gamma _{j}^{(i)}\right), \ \ \ \
-\ell\leq i\leq \ell,\ 1\leq j\leq
k+\max(0,i),\nonumber\\
f^{(s,t)}&=&f(\ve_{s,t}),\ \ \ \ 1\leq
s<t\leq\ell,\label{fnotation}\\
f_i&=&(f_{i,1},f_{i,2},...,f_{i,k+\max(0,i)})\ \ \ \ -\ell\leq i\leq
\ell,\nonumber \\
C(f)&=&(f^{(s,t)})_{1\leq s<t\leq\ell}.\nonumber
\end{eqnarray} When $f\in\Omega$, then $f$ corresponds to the pattern
formed by the rows
$f_{-\ell},f_{-\ell+1},...,f_\ell$ and to the matrix $C=C(f)  \in\Z^{\ell(\ell-1)/2}_{\geq 0}$. Next, we define certain linear
functionals on $\R^{\tG }$ as follows: For $1\leq j\leq\ell$, let
\begin{eqnarray} A_j(f)&=&\sum_{a=1}^{k+j}f_{j,a} -
\sum_{b=1}^{k+j-1}f_{j-1,b},\nonumber \\
B_j(f)&=&  \sum_{a=1}^{k}(f_{-j,a}-f_{-j+1,a}),\label{ABnot}\\
P_j(f)&=&A_j(f)+B_j(f)+\sum_{a<j}f^{(a,j)}+\sum_{b>j}f^{(j,b)}.\nonumber
\end{eqnarray}
We also write
\begin{eqnarray}
A(f)&=&(A_1(f),...,A_\ell(f))\nonumber\\
B(f)&=&(B_1(f),...,B_\ell(f))\label{ABvec}\\
P(f)&=&(P_1(f),...,P_\ell(f)).\nonumber
\end{eqnarray}
For a finite set $X$, $|X|$ shall denotes its cardinality.


\medskip \begin{lemma}\label{come}
For Young diagrams $D$ and $F$ such that $r(D)\leq k$ and $r(F)\leq
k+\ell$, and for $P \in\Z^\ell_{\geq 0}$, let
\[\calc_{F,D,P}=
\{f\in\calc :\ f_{-\ell}=D,f_\ell=F,P(f)=P\}
\] and
\[\Omega_{F,D,P}= \calc_{F,D,P}\cap\Z^{\tG }.\]
Then $\calc_{F,D,P}$ is a polytope and
\[|\Omega_{F,D,P}|=\dim\cale_{F,D,P}.\]
\end{lemma}
\proof For a Young diagram $E$ with $r(E)\leq k$,
$A,B\in\Z^\ell_{\geq 0}$ and $C\in\Z^{\ell(\ell-1)/2}_{\geq 0}$
 such that $S(A,B,C)=P$ (see (\ref{pcon})), let
$\Omega_{F,D,P}^{(E,A,B,C)}$ be the subset of $\Omega_{F,D,P}$
containing all elements $f$ with the properties
\[f_0=E,\ A(f)=A,\ B(f)=B,\ \mbox{and}\ C(f)=C.\]
Then by identifying elements of $\Omega_{F,D,P}^{(E,A,B,C)}$ with patterns of the form \eqref{pattern}, we see that  $|\Omega_{F,D,P}^{(E,A,B,C)} |= K_{F/E,A}K_{D/E,B}$.

\medskip Now $\Omega_{F,D,P}$ is the disjoint union
of the subsets $\Omega_{F,D,P}^{(E,A,B,C)}$. So the cardinality of
$\Omega_{F,D,P}$ is given by the sum of $K_{F/E,A}K_{D/E,B}$ over
all possible $E$, $A$, $B$ and $C$. The number obtained is the
dimension of $\cale_{F,D,P}$ given in Proposition \ref{dimf}. $\Box$

\medskip
Lemma \ref{come} suggests that $\Omega_{F,D,P}$ may be used to index
a basis for $\cale_{F,D,P}$. Since
\[\Omega =\bigcup_{F,D,P}\Omega_{F,D,P},\]
$\Omega$ may be used to index a basis for the algebra $\faol$. We
also observe that if $f\in \Omega_{F,D,P}$ and $g\in
\Omega_{F^\prime,D^\prime,P^\prime}$, then
$f+g\in\Omega_{D+D^\prime,F+F^\prime,P+P^\prime}$. We express this
property as
\begin{equation}\label{grading}
\Omega_{F,D,P}+\Omega_{F^\prime,D^\prime,P^\prime}\subseteq
\Omega_{F+F^\prime, D+D^\prime, P+P^\prime}.
\end{equation}

\medskip

\subsection{The generators of the lattice cone $\Omega$}

In this subsection, we shall use the results in \cite{Ho2} to
specify a set of generators for $\Omega$ and the relations between
the generators.

\medskip
Recall the poset $\tG=\Gamma+F_\ell$ defined in  \eqref{fposet} where $\Gamma=\Gamma(k,\ell)$.
A subset $A$ of a poset  $P$ is said to be \emph{increasing} if
for any $a\in P$,
\begin{equation*}
x\in P \ \mbox{and}\ x \geq a \Longrightarrow x\in A.
\end{equation*}
Let $J^\ast(\tG )$ be the collection of increasing subsets of
$\tG $. For each $A\in J^\ast(\tG )$, let $%
\chi_A:\tG \rightarrow\Z_{\geq 0}$ be the characteristic function of
$A$, that is,
\begin{equation*}
\chi_A(x)=\left\{%
\begin{array}{ll}
1 & x\in A \\
0 & x\not\in A.%
\end{array}%
\right.
\end{equation*}
(If $A=\emptyset$, then $\chi_A(x)=0$ for all $x\in A$.)
Then $\chi_A\in\Omega$. By Theorem 3.3 of \cite{Ho2}, $\Omega $ is
generated by
\begin{equation}\label{calgset}
 \calg=\calg(k,\ell)=\{\chi_A:\ A\in J^\ast(\tG )\},
\end{equation}
and a defining set of relations is given by
\begin{equation*}
\chi_A+\chi_B=\chi_{A\cup B}+\chi_{A\cap B},\ \ \ A,B\in
J^\ast(\tG ) .
\end{equation*}

\medskip
We now specify the increasing sets in $\tG $.
Since $\tG=\Gamma+F_\ell$, every increasing subset $A$ of $\tG$ is of the form $A_1\cup A_2$ where $A_1$ and $A_2$ are increasing subsets of $\Gamma $ and $F_\ell$, respectively.
On the other hand, every subset of $F_\ell$ is increasing.
Thus we only need to determine the increasing subsets of $\Gamma$.
For $0\leq c\leq k$,
$I=\{i_1,...,i_u\}$, $J=\{j_1,...,j_v\}\subseteq \{1,...,\ell\}$
with $|I|=u\leq k-c$, we define a sequence
$(a_{-\ell},...,a_0,...,a_\ell)$ as follows:
 \begin{enumerate}
 \item[(i)] $a_0=c$.
 \item[(ii)] For $0\leq s\leq\ell-1$,
\begin{equation}\label{con1}
a_{-s-1}=\left\{
\begin{array}{ll}
a_{-s}+1 & s+1\in\{i_1,...,i_u\} \\
a_{-s} & s+1\not\in\{i_1,...,i_u\}%
\end{array}
\right.
\end{equation}
and
\begin{equation}\label{con2}
a_{s+1}=\left\{%
\begin{array}{ll}
a_{s}+1 & s+1\in\{j_1,...,j_v\} \\
a_{s} & s+1\not\in\{j_1,...,j_v\}.
\end{array}\right.
\end{equation}
\end{enumerate}
Let \[A_{(c,I,J)}=\{\gamma^{(i)}_j: 1\leq j\leq a_i,\ -\ell\leq
i\leq \ell\}.\]
(So $A_{0,\emptyset,\emptyset}=\emptyset$.)
Then one verifies that $A_{(c,I,J)}$  is an increasing subset of $\Gamma$. In fact, the sets $A_{(c,I,J)}$ exhaust all the increasing subsets of $\Gamma$.

\medskip
\begin{lemma}\label{ics}  Every increasing subset of $\Gamma$ is of the form $A_{(c,I,J)}$. Consequently,
\[ J^\ast(\tG )=
\{A_{(c,I,J)}\cup Z:\ 0\leq c\leq k,\ I,J\subseteq\{1,...,\ell\},
|I|\leq c-k, Z\subseteq F_\ell  \}. \]
\end{lemma}

\proof Let $A$ be an increasing subset of $\Gamma$. Then for $-\ell\leq
j\leq \ell$,
\[A\cap\{\gamma^{(j)}_i: 1\leq i\leq k+\max(0,j)\}=\{\gamma^{(j)}_i:\ 1\leq i\leq a_j\}\]
for some $a_j\geq 0$. ($a_j=0$ means the intersection is empty.)

\medskip Let $0\leq s\leq \ell-1$. Then
$\gamma^{(s)}_{a_s}\in A$ and $\gamma^{(s)}_{a_s+1}\not\in A$. Since
$A$ is increasing and  $\gamma^{(s+1)}_{a_s}\geq
\gamma^{(s)}_{a_s}$, $\gamma^{(s+1)}_{a_s}\in A$. Similarly, since
$\gamma^{(s)}_{a_s+1}\geq \gamma^{(s+1)}_{a_s+2}$, we must have
$\gamma^{(s+1)}_{a_s+2}\not\in A$. So $a_{s+1}=a_s$ or
$a_{s+1}=a_s+1$. Similarly, $a_{-s-1}=a_{-s}$ or
$a_{-s-1}=a_{-s}+1$. We now let $c=a_0$ and define $I$ and $J$ by
equations \eqref{con1} and \eqref{con2}. Then $A=A_{(c,I,J)}$.
 $\Box$

\subsection{Standard monomial basis for $\faol$}\label{smb}
Recall that $\faol=\left(\calp_{n,k,\ell}\right)^{U_{n} \times
U_{k}}$ where
\[\calp_{n,k,\ell}= \pmnk\otimes \left(\bigotimes_{i=1}^\ell \calp(\C
^n_i)\right)\otimes \left(\bigotimes_{j=1}^\ell \calp(\C^k_j)\right)\otimes%
\left(\bigotimes_{1\leq s<t\leq \ell} \calp(\C_{s,t})\right).\]
 We denote the standard coordinates on
$\Mnk$, $\C^n_j$, $\C^k_j$ and $\C_{s,t}$ as follows:
\begin{enumerate}
\item[(i)] $(x_{ij})_{1\leq i\leq n,1\leq j\leq k}\in\Mnk$.
\item[(ii)] $(y_{ij})_{1\leq i\leq n}\in \C^n_j$, $1\leq j\leq \ell$.
\item[(iii)] $(r_{i,k+j})_{1\leq i\leq k}\in\C^k_j$, $1\leq j\leq \ell$.
\item[(iv)] $(r_{k+s,k+t})\in \C_{s,t}$, $1\leq s<t\leq\ell$.
\end{enumerate}
So $\calp_{n,k,\ell}$ can be regarded as a polynomial algebra on
these variables.

\medskip For later use, we define a monomial ordering
on  $\calp_{n,k,\ell}$ as follows: it is the graded lexicographic
 order (\cite{CLO}) with respect to
\begin{eqnarray*}
x_{11} &>&x_{21}>\cdots >x_{n1}>x_{12}>x_{22}>\cdots >x_{nk} \\
&>&y_{11}>y_{21}>\cdots >y_{n1}>y_{21}>\cdots >y_{n\ell } \\
&>&r_{1,k+1}>r_{2,k+1}>\cdots >r_{k,k+1}>r_{1,k+2}>\cdots
>r_{k,k+\ell} \\
&>&r_{k+1,k+2}>r_{k+1,k+3}>r_{k+2,k+3}>\cdots
>r_{k+\ell -1,k+\ell. }
\end{eqnarray*}
For each polynomial $f$ in $\calp_{n,k,\ell}$, we write the leading
monomial of $f$ with respect to the above monomial ordering as
$\LT(f)$.

\medskip

\noindent \textbf{Definition.} For $0\leq c\leq k$,
$I=\{i_{1}<i_{2}<\cdots <i_{u}\}, \ J=\{j_{1}<j_{2}<\cdots
<j_{v}\}\subseteq \{1,...,\ell \}$ with $|I|=u\leq k-c$,  let
\begin{equation*}
\eta _{(c,I,J)}=\left\vert
\begin{array}{ccccccc}
x_{11} & x_{12} & \cdots  & x_{1,c+u} & y_{1,j_{1}} & \cdots  &
y_{1,j_{v}}
\\
x_{21} & x_{22} & \cdots  & x_{2,c+u} & y_{2,j_{1}} & \cdots  &
y_{2,j_{v}}
\\
\vdots  & \vdots  &  & \vdots  & \vdots  &  & \vdots  \\
x_{c+v,1} & x_{c+v,2} & \cdots  & x_{c+v,c+u} & y_{c+v,j_{1}} &
\cdots  &
y_{c+v,j_{v}} \\
r_{1,k+i_{1}} & r_{2,k+i_{1}} & \cdots  & r_{c+u,k+i_{1}} & 0 & \cdots  & 0 \\
\vdots  & \vdots  &  & \vdots  & \vdots  &  & \vdots  \\
r_{1,k+i_{u}} & r_{2,k+i_{u}} & \cdots  & r_{c+u,k+i_{u}} & 0 & \cdots  & 0%
\end{array}%
\right\vert .
\end{equation*}%
(So $\eta_{(0,\emptyset,\emptyset)}=1$.)

\medskip For each $g\in\calg$, we define the polynomial
$\eta_g\in\faol$ as follows:
\begin{enumerate}
\item[(i)] If $g=\chi_{A_{(c,I,J)}}$, then $\eta_g=\eta_{(c,I,J)}$.
\item[(ii)] If $g=\chi_{\{\ve_{s,t}\}}$, then $\eta_g=r_{k+s,k+t}$.
\item[(iii)] If $g=\chi_B$ where $B=A_{(c,I,J)}\cup Z$ and
$Z\subseteq \Gamma_\ell$, then
\[\eta_g=\eta_{(c,I,J)}\cdot\prod_{\ve_{s,t}\in Z}r_{k+s,k+t}.\]
\end{enumerate}

\medskip More generally, we can extend the map $g\rightarrow
\eta_g$ on $\calg $ to $\Omega $. This is done in the following way.
Recall from \cite{Ho2} and \cite{Hi} that each element $g$ has a
unique standard expression as a sum
\begin{equation}\label{eqg}
g=\sum_{j=1}^{N}c_j\chi_{A_j}
\end{equation}
where $N$ is given in equation \eqref{von}, $c_j \geq 0$ and
$\emptyset \subset A_1 \subset A_2 \subset \cdots \subset A_N=\tG $
is a maximal chain in $J^\ast(\tG )$. Then we define
\begin{equation}\label{eqgamma}
\eta_g=\prod_{j=1}^{N}\eta_{\chi_{A_j}}^{c_j}.
\end{equation}
Let
\[\calb=\{\eta_g:\ g\in\Omega \}.\]
Note that $\calb$ can be characterized as follows. Let
\[S=S(k,\ell)=\{\eta_g:\ g\in\calg \}.\] Then the partial
ordering on $J^\ast(\tG )  $ induces the following partial ordering
on $S $: if $g_1,g_2\in S $, then $g_1=\chi_{A_1}$ and
$g_2=\chi_{A_2}$ for some $A_1,A_2\in {J^\ast(\tG )}$. We define
\[g_1\leq g_2\Longleftrightarrow A_1\subseteq A_2.\]
Since $J^\ast(\tG )$ is a distributive lattice, so is $S $. We now
observe that each $\eta_g\in\calb$ is a monomial in elements of a
maximal chains in $S$, that is, it is a {\em standard monomial}. We
shall show that $\calb$ is a basis for $\faol$, so that the algebra
$\faol$ has a standard monomial theory for $S$.
  In the appendix, we will give the Hasse diagram of the poset
$J^\ast(\tG )$ for $\ell=1$.

\medskip
\begin{lemma} The leading monomial of $\eta_{(c,I,J)}$ is given by
\[\LT(\eta_{(c,I,J)})=\left(\prod_{a=1}^cx_{aa}\right)
\left(\prod_{a=1}^vy_{c+a,j_a}\right)
\left(\prod_{a=1}^ur_{c+a,k+i_a}\right).
\]
\end{lemma}
\proof By observation. $\Box$

\medskip\begin{prop}\label{LT1}
For $g\in \Omega $, with notation given in (\ref{fnotation}),
\begin{equation}\label{ltg}
\LT(\eta_g)=\left(\prod_{u=1}^kx_{u,u}^{g_{0,i}}\right)
 \left(\prod_{
\substack{ 1\leq a\leq k+b \\ 1\leq b\leq \ell
}}y_{a,b}^{g_{b,a}-g_{b-1,a}} \right) \left( \prod_{\substack{ 1\leq
i\leq k \\ 1\leq j\leq \ell }}r_{i,k+j}^{g_{-j,i}-g_{-j+1,i}}
\right) \left( \prod_{1\leq s<t\leq \ell}r_{k+s,k+t}^{g^{(s,t)}}
\right).
 \end{equation}
\end{prop}
\proof For each $g\in \Omega$, let $m(g)$ be the monomial given in
the right hand side of  \eqref{ltg}. This defines a map
$m$ from $\Omega$ to the semigroup of all monomials in $\calp_{n,k,\ell}$. It is easy to check that $m$
is a semigroup homomorphism, and $m(g)=\LT(\eta_g)$ for all
$g\in\calg$. In fact, $m$ is an semigroup isomorphism onto its
image.

\medskip We now let $g\in\Omega$ and
assume that its standard form is given in equation \eqref{eqg}. Then
$\eta_g$ is given in equation \eqref{eqgamma}. Since $m$ is a
semigroup homomorphism,
\[\LT(\eta_g)=\prod_{j=1}^{N}\LT(\eta_{\chi_{A_j}})^{c_j}
=\prod_{j=1}^{N}m\left(\chi_{A_j}\right)^{c_j}
=m\left(\sum_{j=1}^{N}c_j\chi_{A_j} \right)=m(g).\ \Box\]

\bigskip\begin{cor}\label{sgis} The set
\[\LT(\calb)=\{\LT(\eta_g):\ g\in \Omega\}\]
of monomials forms an affine semigroup isomorphic to $\Omega$.
\end{cor}
\proof This is because the map $m(g)=\LT(\eta_g)$ defined in the
proof Proposition \ref{LT1} is an semigroup isomorphism from $\Omega
$ onto $\LT(\calb)$. $\Box$

\bigskip\begin{cor} The set
$\calb$  is linearly independent.
\end{cor}
\proof This is because
the elements in $\calb$ have distinct leading monomials. $\Box$

\bigskip\begin{prop}\label{tnk1}
Let
\[\calb_{F,D,P}=\{\eta_g:\ g\in\Omega_{F,D,P}\}.\]
Then
  $\calb_{F,D,P} $ is a basis
for the homogeneous component $\cale_{F,D,P}$ of $\faol$.
 \end{prop}

\proof We claim that if $g\in \Omega_{F,D,P}$, then
$\eta_g\in\cale_{F,D,P}$. This is true if $g\in\calg$. The general
case follows from this and    \eqref{grading}. It follows that
$\calb_{F,D,P}\subseteq \cale_{F,D,P}$. Moreover, $\calb_{F,D,P}$ is
linearly independent and the number of vectors it contains coincides
with the dimension of $\cale_{F,D,P}$. So $\calb_{F,D,P}$ is a basis
for $\cale_{F,D,P}$. $\Box$

\medskip\begin{thm} \label{smb} The algebra $\faol$ has a standard monomial
theory for $S$, that is,  the monomials in elements of the maximal
chains in $S$ form a basis for $\faol$.
\end{thm}

\proof The set $\calb$ is precisely the set of
  all monomials in elements of
the maximal chains in $S$. It forms a basis  for $\faol$ follows
from $\calb=\bigcup_{F,D,P}\calb_{F,D,P}$,   \eqref{homcom} and
Proposition \ref{tnk1}.
   $\Box$

\subsection{Toric degeneration of $\faol$}
We first recall the notion of a {\em SAGBI basis} (\cite{RS, Stu}).
Let $R=k[x_1,...,x_n]$ be a polynomial algebra over a
field $k$ and $L$ a subalgebra of  $R$. Assume that $R$ is given a
monomial ordering. For each $a\in R$, let $\LT(a)$ be the leading
monomial of $a$ with respect to this monomial ordering. The {\em
initial algebra of } $L$, denoted by $\LT(L)$, is the subalgebra of
$R$ generated by the set $\{\LT(a):\ a\in L\}$ of leading monomials.
A subset $F$ of $L$ is called a {\em SAGBI basis} for $L$ if the set
$\{\LT(a):\ a\in F\}$ generates the initial algebra $\LT(L)$. We
will need the following result.

\medskip
\begin{prop} \label{chv} {\cite{CHV}} Let $L$ be a
subalgebra of a polynomial algebra $R=k[x_1, ..., x_n]$ over a field
$k$.  If the initial algebra $\LT(L)$ of $L$ with respect to some
monomial ordering is finitely generated, then there exists a flat
$1$-parameter family of $k$-algebras with general fibre $L$ and
special fibre $\LT(L)$.
\end{prop}

\medskip Let $\C[\Omega ]$ be the semigroup algebra (\cite{BH}) on
$\Omega $, that is, $\C[\Omega]$ is a complex algebra with basis
$\{X^f:\ f\in\Omega\}$ and it satisfies the multiplication law
\[X^{f_1}X^{f_2}=X^{f_1+f_2},\ \ \ f_1,f_2\in\Omega.\]
On the other hand, the distributive lattice $J^\ast(\tG)$ gives rise
to the Hibi algebra (\cite{Hi, Ho2}) generated by
$J^\ast(\tG)$ and with relations
\[A_1\cdot A_2=(A_1\cup A_2)\cdot(A_1\cap A_2).\]
It is easy to see that this algebra and the semigroup algebra
$\C[\Omega]$ are isomorphic (\cite{Ho2}).

\bigskip\begin{lemma}\label{sgis} The initial algebra
$\LT(\faol)$ of $\faol$ is isomorphic to the semigroup algebra
$\C[\Omega]$.
\end{lemma}
\proof Let $f\in\faol$. By Theorem \ref{smb},
\[f=c_1\eta_{g_1}+c_2\eta_{g_2}+\cdots+c_r\eta_{g_r}\]
for some $c_1,...,c_r\in\C$ and $g_1,...,g_r\in\Omega$. Since the
leading monomials $\LT(\eta_{g_1}), ..., $ $\LT(\eta_{g_r})$ are
distinct, $\LT(f)=\LT(\eta_{g_j})$ for some $1\leq j\leq r$. This
shows that $\LT(f)\in\LT(\calb)$. It follows that the initial
algebra $\LT(\faol)$ is generated by $\LT(\calb)$. Now by Lemma
\ref{sgis},  the semigroups $\LT(\calb)$ and $\Omega$   are
isomorphic. So the algebras $\LT(\faol)$ and $\C[\Omega]$ are also
isomorphic. $\Box$

\medskip

\medskip\begin{thm} There exists a flat one-parameter family of
complex algebras with general fibre $\faol$ and special fibre
$\C[\Omega]$.
\end{thm}

\proof Since $\calg$ generates the semigroup $\Omega$,
  the set $\{\LT(\eta_f):\ \eta_f\in S\}$ generates the
initial algebra $\LT(\faol)$. So $S$ is a finite SAGBI
basis for $\faol$. The theorem now follows from Proposition
\ref{chv} and Lemma \ref{sgis}. $\Box$

\section{Pieri algebras for symplectic groups}
There is a parallel construction of  Pieri algebras for the complex
symplectic group $\Sptn$ which we will briefly explain in this
section. This algebra encodes information on how a tensor product of
$\Sptn$ representations of the form
\[\tau^D_{2n}\otimes\tau^{(p_1)}_{2n}
\otimes\tau^{(p_2)}_{2n}\otimes\cdots \otimes\tau^{(p_\ell)}_{2n}\]
decomposes, where $r(D)\leq k$, $p_1,...,p_\ell\geq 0$ and $k+\ell\leq n$ (see
Section \ref{spr} for notation).

\medskip
Let the groups $\Sptn$ and $\Glk$ act on the algebra $\pmtnk$
of polynomial functions on $\Mtnk$ by
\begin{equation}\label{spact}
 [(g,h).f](T)=f(g^tTh),\ \ \ g\in\Sptn,\ h\in\Glk,\ f\in\pmtnk,\ T\in\Mtnk.
 \end{equation}
Let $I_{2n,k}$ be the ideal of $\pmtnk$ generated by all $\Sptn$ invariants of positive degree, and form the quotient algebra $\pmtnk/I_{2n,k}$. Then the action \eqref{spact} induces an action of $\Sptn\times \Glk$ on $\pmtnk/I_{2n,k}$, and
\[\pmtnk/I_{2n,k}\cong\sum_{r(D)\leq k}\tau^D_{2n}\otimes\rho^D_k.\]
Let $\left(\pmtnk/I_{2n,k}\right)^{U_k}$ be the subalgebra of $U_k$ invariants in $\pmtnk/I_{2n,k}$. It is a module for $\Sptn\times A_k$.

\medskip Next, for each $1\leq j\leq\ell$, let
$\C^{2n}_j$ (resp. $\GL^{(j)}_1$) be a copy of $\C^{2n}$ (resp. $\GL_1$) and consider the algebra
$\calp(\C^{2n}_j)$. Since $\C^{2n}_j\cong\rM_{2n,1}$, by the $(\GL_{2n},\GL^{(j)}_1)$-duality,
\[\calp(\C^{2n}_j)\cong\sum_{p_j\geq 0}\rho^{(p_j)}_{2n}\otimes\rho^{(p_j)}_1\]
as a $\GL_{2n}\times\GL^{(j)}_1$ module. We now restrict the action of $\GL_{2n}$ to $\Sptn$. But for each $p_j\geq 0$, $\rho^{(p_j)}_{2n}$ remains irreducible under the action by $\Sptn$, so that $\rho^{(p_j)}_{2n}=\tau^{(p_j)}_{2n}$. Thus
\[\calp(\C^{2n}_j)\cong\sum_{p_j\geq 0}\tau^{(p_j)}_{2n}\otimes\rho^{(p_j)}_{1,j}.\]
as a $\Sptn\times\GL^{(j)}_1$ module.

\medskip We now let
\begin{eqnarray*}
\asp&:=&\left\{ \left( \pmtnk/I_{2n,k} \right)^{U_k}\otimes
 \calp(\C^{2n}_1) \otimes\cdots\otimes
  \calp(\C^{2n}_\ell)
 \right\}^{U_{\Sptn}}\\
 &=& \left\{
\left(    \pmtnk/I_{2n,k} \right)\otimes
 \calp(\C^{2n}_1)\otimes\cdots\otimes
  \calp(\C^{2n}_\ell)
 \right\}^{U_{\Sptn}\times U_k}.
\end{eqnarray*}
This algebra  is a module for $A_{\Sptn}\times A_k\times \GL^{(1)}_1\times\cdots\times \GL_1^{(\ell)}$
and can be decomposed as
\[ \asp\cong
\sum_{\stackrel{\scriptstyle r(D)\leq k}{p_1,...,p_\ell\geq
0}}\left(\tau^D_{n}\otimes \tau^{(p_1)}_{n}\otimes\cdots\otimes
\tau^{(p_\ell)}_{n}\right)^{U_{\Sptn}}\otimes \left(
\rho^D_k\right)^{U_k}\otimes \rho^{(p_1)}_{1,1}\otimes\cdots\otimes
\rho^{(p_\ell)}_{1,\ell}.\]
For Young diagrams $D$ and $F$ with $r(D)\leq k$ and $r(F)\leq k+\ell$, and $P=(p_1,...,p_\ell)\in\Z^\ell_{\geq 0}$,
the $\chi^F_n\times\psi^D_k\times\psi^{(p_1)}_1\times\cdots\times
\psi^{(p_\ell)}_1$-eigenspace
in this algebra can be identified with the space of all
$\Sptn$ highest weight vectors of weight $\chi^F_n$ in
$\tau^D_{2n}\otimes \tau^{(p_1)}_{2n}\otimes\cdots\otimes
\tau^{(p_\ell)}_{2n}$. Thus the multiplicity of $\tau^F_{2n}$ in this tensor product coincides with the dimension of the eigenspace.
  In view of this property, we will call $\asp$ an {\em   $\Sptn$ Pieri
algebra.}

\medskip Using similar arguments as in Proposition \ref{ispa},
  we can show that
the algebras $\asp$ and $\mathfrak{A}_{2n,k,\ell}$ are isomorphic. Consequently, one can deduce the structure of $\asp$ from the results of Section \ref{stfaol}.
In particular, $\asp$ has a standard monomial
basis and has a flat deformation to a Hibi algebra.

\section{Appendix}
We show the poset structure of the generating set $S(k,\ell)$ of the algebra $\faol$ for the case $\ell=1$. Note that if $\ell=1$, then we have $\tG(k,\ell)=\Gamma(k,\ell)$. Therefore, from Lemma \ref{ics}, the poset structure of $S(k,\ell)$ can be obtained by considering order increasing subsets of $\Gamma(k,1)$ drawn in \eqref{Gamma ell=1}.

Following the notation used in \eqref{con1} and \eqref{con2}, we let $(a_{-1},a_0,a_1)$ denote the order preserving characteristic function $\chi$ over $\Gamma(k,1)$ such that the integers $a_{-1}$, $a_0$, and $a_1$ are equal to the cardinalities of the intersections between $\chi^{-1}(1)$ and the top row, the middle row, and the bottom row of $\Gamma(k,1)$ respectively.

\newpage

{\tiny
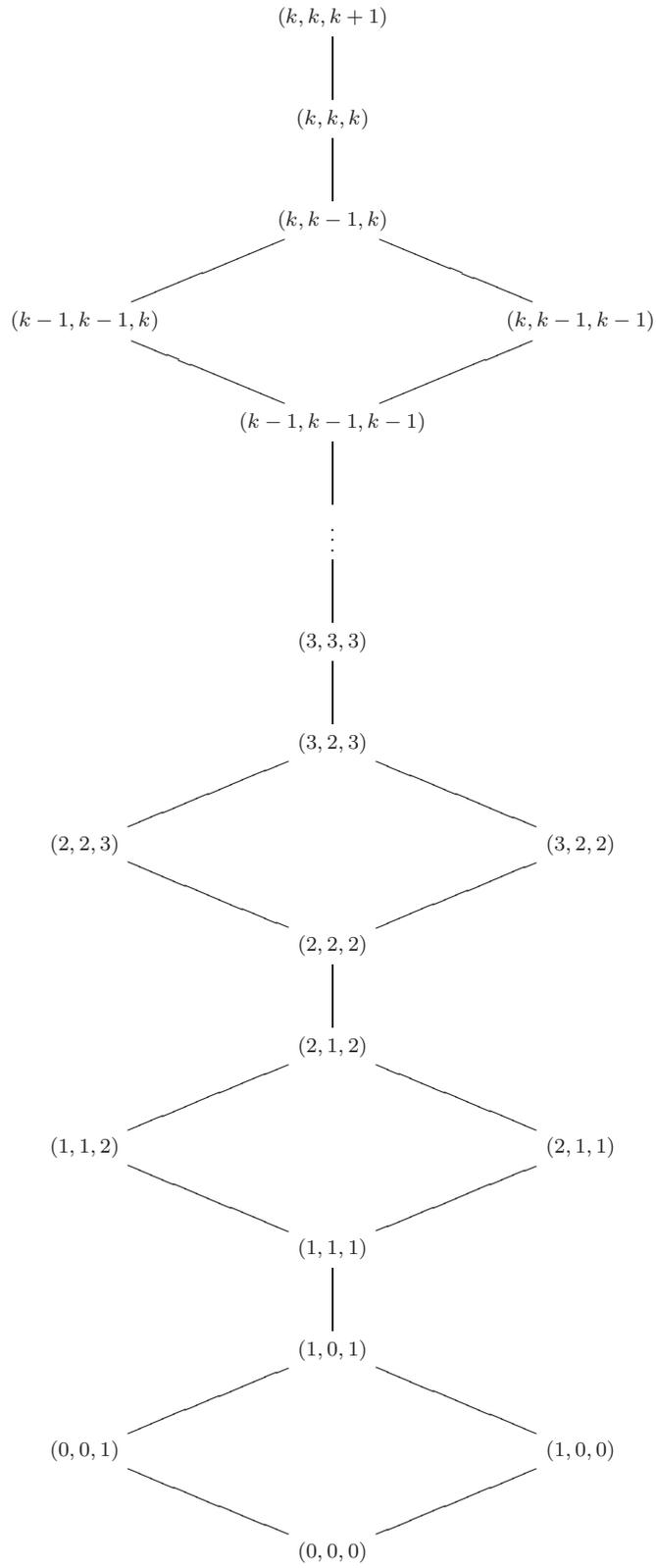
\begin{equation*}
\xymatrix{%
                       & (k,k,k+1) \ar@{-}[d]           &  \\
                       & (k,k,k)   \ar@{-}[d]           &  \\
                       & (k,k-1,k)                      &  \\
(k-1,k-1,k) \ar@{-}[ur]&                                & (k,k-1,k-1) \ar@{-}[ul]\\
                       & \ar@{-}[ul] (k-1,k-1,k-1) \ar@{-}[d] \ar@{-}[ur]&  \\
                       & \vdots   \ar@{-}[d]               & \\
                       & (3,3,3)   \ar@{-}[d]           &  \\
                       & (3,2,3)                      &  \\
(2,2,3) \ar@{-}[ur]&                                & (3,2,2) \ar@{-}[ul]\\
                       & \ar@{-}[ul] (2,2,2) \ar@{-}[d] \ar@{-}[ur]&  \\
                       & (2,1,2)                      &  \\
(1,1,2) \ar@{-}[ur]&                                & (2,1,1) \ar@{-}[ul]\\
                       & \ar@{-}[ul] (1,1,1) \ar@{-}[d] \ar@{-}[ur]&  \\
                       & (1,0,1)                        &  \\
(0,0,1) \ar@{-}[ur]&                                    & (1,0,0) \ar@{-}[ul]\\
                       & \ar@{-}[ul] (0,0,0) \ar@{-}[ur] &  }
\end{equation*}}


\begin{thebibliography}{LLL99}

\bibitem[AB]{AB} V. Alexeev and M. Brion, \textit{Toric degenerations of
spherical varieties}, Selecta Math. (N.S.) 10 (2004), 453 -- 478.

\bibitem[BH]{BH} W. Bruns and J. Herzog, \textit{Cohen-Macaulay
Rings}, Cambridge Univ. Press, 1998, 453 pp.

\bibitem[Ca]{Ca} P. Caldero, \textit{Toric degenerations of Schubert varieties},
Transform. Groups 7 (2002), no. 1, 51 -- 60.

\bibitem[CHV]{CHV} A. Conca, J. Herzog and G. Valla,
\textit{SAGBI basis with applications to blow-up algebras}, J. Reine
Angew. Math. 474 (1996), 113 -- 138.

\bibitem[CLO]{CLO} D. Cox, J. Little and D. O'Shea,
\textit{Ideals, Varieties, and Algorithms. An Introduction to
Computational Algebraic Geometry and Commutative Algebra}, Second
edition, Undergraduate Texts in Mathematics, Springer-Verlag, New
York, 1997.



\bibitem[Fu]{Fu2} W. Fulton, \textit{Young Tableaux},
Cambridge University Press, Cambridge, UK, 1997.

\bibitem[GL]{GL} N. Gonciulea and V. Lakshmibai, Degenerations of flag and
Schubert varieties to toric varieties. Transform. Groups 1 (1996), no. 3,
215--248.

\bibitem[GW]{GW} R. Goodman and N. Wallach,
\textit{Representations and Invariants of the Classical Groups},
Cambridge Univ. Press, 1998.

\bibitem[Hi]{Hi} T. Hibi, {\it Distributive lattices, affine semigroup rings and
algebras with straightening laws}, in: Commutative Algebra and
Combinatorics, in: Adv. Stud. Pure Math., vol. 11, 1987, pp. 93–109.

\bibitem[Ho1]{Ho1} R. Howe,
 \textit{ Perspectives on invariant theory},
The Schur Lectures, I. Piatetski-Shapiro and S. Gelbart (eds.),
Israel Mathematical Conference Proceedings, 1995, 1 -- 182.

\bibitem[Ho2]{Ho2} R. Howe, \textit{
Weyl Chambers and Standard Monomial Theory for Poset Lattice Cones},
Q. J. Pure Appl. Math. 1 (2005) 227–239.

\bibitem[HJL${}^+$]{HJLTW} R. Howe, S. Jackson, S. T. Lee, E-C Tan
and J. Willenbring, \textit{Toric degeneration of branching
algebras}, Adv. in Math. 220 (2009), 1809--1841.


\bibitem[HL1]{HL1} R. Howe and S. T. Lee, \textit{Bases for
some reciprocity algebras I}, Trans. Amer. Math. Soc. 359 (2007),
4359-4387.

\bibitem[HL2]{HL2} R. Howe and S. T. Lee, \textit{Bases for
some reciprocity algebras II}, Adv. in Math. 206 (2006), 145-210.

\bibitem[HL3]{HL3} R. Howe and S. T. Lee, \textit{Bases for
some reciprocity algebras III}, Compositio Math. 142 (2006),
1594-1614.

\bibitem[HTW1]{HTW1} R. Howe, E-C. Tan and J. Willenbring,
\textit{Reciprocity algebras and branching for classical symmetric
pairs}, in  ``Groups and Analysis - the Legacy of Hermann Weyl,"
London Mathematical Society Lecture Note Series No. 354, Cambridge
University Press, 2008.

\bibitem[HTW2]{HTW3} R. Howe, E-C. Tan and J. Willenbring,
\textit{A basis for $GL_n$ tensor product algebra}, Adv. in Math.
196 (2005), 531 -- 564.



\bibitem[Ki1]{Ki1} S. Kim,
{\it Standard monomial theory for flag algebras of $GL(n)$ and
$Sp(2n)$,} J. Algebra 320 (2008), 534 -- 568.

\bibitem[Ki2]{Ki2} S. Kim,
{\it Standard monomial bases and degenerations of $SO(m)$
representations,} J. Algebra 322(11) (2009), 3896 -- 3911

\bibitem[KM]{KM} M. Kogan and E. Miller, Toric degeneration of Schubert
varieties and Gelfand-Tsetlin polytopes. Adv. Math. 193 (2005), no. 1, 1 -- 17.

\bibitem[La]{La} V. Lakshmibai, The development of standard monomial theory.
II. A tribute to C. S. Seshadri (Chennai, 2002), 283--309, Trends Math.,
Birkh\"{a}user, Basel, 2003.

\bibitem[Mu]{Mu} C. Musili, The development of standard monomial theory. I. A
tribute to C. S. Seshadri (Chennai, 2002), 385--420, Trends Math., Birkh\"{a}%
user, Basel, 2003.

\bibitem[RS]{RS} L. Rubbiano and M. Sweedler, \textit{Subalgebra bases},
Proceedings Salvador 1988 Eds. W. Bruns, A. Simis, Lecture Notes
Math 1430 (1990), 61 -- 87.

\bibitem[Sta]{Sta} R. Stanley, \textit{Enumerative combinatorics,
Volume 2}, Cambridge University Press, 1999.

\bibitem[Stu]{Stu} B. Sturmfels,
\textit{Grobner Beses and Convex Polytopes}, Univ. Lecture Series,
Vol. 8, Amer. Math. Soc., Providence, RI, 1996.

\bibitem[Wy]{Wy} H. Weyl, \textit{The Classical Groups}, Princeton University
Press, Princeton, N.J., 1946.

\end{thebibliography}
\end{document}